\documentstyle{amsppt} \magnification=1200 \NoBlackBoxes

\topmatter \title RANDOM POLYTOPES AND AFFINE SURFACE AREA \endtitle
\author Carsten Sch\"utt \endauthor \address Oklahoma State University,
Department of Mathematics, Stillwater, Oklahoma 74078 \endaddress
\address Christian-Albrechts Universit\"at, Mathematisches Seminar,
24118 Kiel, Germany \endaddress \thanks  Supported by NSF-grant
DMS-9301506 \endthanks \subjclass 52A22 \endsubjclass \abstract Let K
be a convex body in $R^d$. A random polytope is the convex hull
$[x_1,...,x_n]$ of finitely many points chosen at random in K. $\Bbb
E(K,n)$ is the expectation of the volume of a random polytope of n
randomly chosen points. I. B\'ar\'any showed that we have for convex
bodies with $C^3$ boundary and everywhere positive curvature

$$ c(d)\lim_{n \to \infty} \frac {vol_d(K)-\Bbb
E(K,n)}{(\frac{vol_d(K)}{n})^{\frac{2}{d+1}}} =\int_{\partial K}
\kappa(x)^{\frac{1}{d+1}}d\mu(x)                 $$

where $\kappa(x)$ denotes the Gau\ss-Kronecker curvature. We show that
the same formula holds for all convex bodies if $\kappa(x)$ denotes the
generalized Gau\ss-Kronecker curvature.  \endabstract \endtopmatter
\document \break

\specialhead 1. Introduction \endspecialhead

Let K be a convex body in $\Bbb R^d$. A random polytope in K is the
convex hull of finitely many points in K that are chosen at random with
respect to a probability measure on K. Here we consider the normalized
Lebesgue measure on K.  For a fixed number n of points we are
interested in the expectation of the volume of that part of K that is
not contained in the convex hull $[x_1,....., x_n]$  of the chosen
points. We denote

$$ \Bbb E(K,n)= \int_{K \times \cdots \times K}
vol_d([x_1,...,x_n])d\Bbb P(x_1,...x_n)   $$

where $\Bbb P$ is the n-fold product of the normalized Lebesgue measure
on K.  We are interested in the asymptotic behavior of

$$ vol_d(K)-\Bbb E(K,n)= \int_{K \times \cdots \times K} vol_d(K
\setminus [x_1,....,x_n]) d\Bbb P(x_1,...,x_n)     $$

In [R-S1, R-S2] the asymptotic behavior of this expression has been
determined for polygons and smooth convex bodies in $R^2$.  \vskip 1cm

\proclaim{\smc Theorem 1} Let K be a convex body in $\Bbb R^d$. Then we
have

$$ c(d)\lim_{n \to \infty} \frac {vol_d(K)-\Bbb
E(K,n)}{(\frac{vol_d(K)}{n})^{\frac{2}{d+1}}} =\int_{\partial K}
\kappa(x)^{\frac{1}{d+1}}d\mu(x)                 $$

where $\kappa (x)$ is the generalized Gau\ss-Kronecker curvature and

$$ c(d)=2(\frac{vol_{d-1}(B_2^{d-1})}{d+1})^{\frac {2}{d+1}}
\frac{(d+3)(d+1)!}
    {(d^2+d+2)(d^2+1)\Gamma(\frac{d^2+1}{d+1})}
    $$

\endproclaim \vskip 1cm

This problem was posed by Schneider and Wieacker [Schn-W] and Schneider
[Schn]. It has been solved by B\'ar\'any [B] for convex bodies with
$C^3$ boundary and everywhere positive curvature. Our result holds for
arbitrary convex bodies. \par The main ingredients of the proof are
taken from [B-L],[B], and [Sch\"u-W 1].  \vskip 1cm

We introduce the notion of generalized Gau\ss-Kronecker curvature. A
convex function $f:X \rightarrow \Bbb R, X \subseteq \Bbb R^d$ is
called twice differentiable at $x_0$ in a generalized sense if there
are a linear map $d^{2}f(x_0) \in L(\Bbb R^d)$ and a neighborhood
$U(x_0)$ so that we have for all $x \in U(x_0)$ and all
subdifferentials df(x)

$$   \parallel df(x)-df(x_0)-(d^{2}f(x_0))(x-x_0) \parallel
     \leq \Theta(\parallel x-x_0 \parallel ) \parallel x-x_0
     \parallel      $$

where $\lim_{t \to 0} \Theta(t)=\Theta(0)=0$ and where $\Theta$ is a
montone function. $d^{2}f(x_0)$ is symmetric and positive semidefinite.
If f(0)=0 and df(0)=0 then the ellipsoid or elliptical cylinder

$$ x^{t}d^{2}f(0)x=1    $$

is called the indicatrix of Dupin at 0. The general case is reduced to
the case f(0)=0 and df(0)=0 by an affine transform. The eigenvalues of
$d^{2}f(0)$ are called the pricipal curvatures and their product the
Gau\ss-Kronecker curvature $\kappa(0)$. Aleksandrov [A, Ba] proved that
a convex surface is almost everywhere differentiable in the generalized
sense. As surface measure on $\partial K$ we take the restriction of
the (d-1)-dimensional Hausdorff measure to $\partial K$. For $x \in
\partial K$ the normal at x to $\partial K$ is denoted by N(x). N(x) is
almost everywhere unique.  We denote

$$ K_t=\{x \in K| vol_d((-x+K) \cap (x-K)) \geq t
\}                      $$

for $t \in [0,T]$ with

$$ T= \max_{y \in K} vol_d((-y+K) \cap
(y-K))                               $$

$K_t$ is a convex body and was studied and used in [St, F-R, K, Schm1]
and was called convolution body in [K,Schm1]. It was shown in [St, F-R]
that $K_T$ consists of one point only.  Therefore we may also interpret
$K_T$ - in abuse of notation - as a point. For a given $x \in \partial
K$ there is a unique point $x_{t} \in \partial K_{t}$ that is the
intersection of the interval $[K_T,x]$ and $\partial K_{t}$. \par
$P_{\xi}$ denotes the orthogonal projection onto the hyperplane
orthogonal to $\xi$ and passing through the origin. $B_2^{d}(x,r)$ is
the Euclidean ball in $R^d$ with center x and radius r.  $B_2^{d}$ is
the ball with center 0 and radius 1. $H(x,\xi)$ denotes the hyperplane
through x and orthogonal to $\xi$. For a given hyperplane H the closed
halfspaces are denoted by $H^{+}$ and $H^{-}$. Usually $H^{+}$ is the
halfspace containing $K_T$ if we consider a convex body K.
 \par \vskip 1.5cm

\specialhead 2. Outline Of The Proof  \endspecialhead

We outline the proof of Theorem 1. We have that

$$ vol_d(K)-\Bbb E(K,n)= \int_{K} \Bbb P\{(x_1,...,x_n)|x \notin
[x_1,...,x_n] \} dx     $$

$$  =- \int_{0}^{T} \frac{d}{dt} \int_{K_t} \Bbb P\{x_1,...,x_n)|x
\notin [x_1,...,x_n] \} dxdt     $$

The derivative can be computed and we get

$$ \int_0^T \int_{\partial K_t} \frac {\Bbb P\{(x_1,...,x_n) \mid x
\notin [x_1,...,x_n]\}}{ vol_{d-1}(P_{N(x)} ((-x+K) \cap (x-K)))}
d\mu_t(x)dt           $$

where $\mu_t$ is the surface measure on $\partial K_t$. We pass to an
integral on $\partial K$ instead of $\partial K_t$.

$$ \int_{\partial K} \int_0^T \frac{\Bbb P \{(x_1,...., x_n) \mid x_t
\notin [x_1,....,x_n] \} } {vol_{d-1}(P_{N(x_t)} ((-x_t+K) \cap
(x_{t}-K)))}  \frac {\Vert x_t \Vert^d}{\Vert x \Vert ^d}  \frac
{<x,N(x)>}{<x_t,N(x_t)>} dt d\mu(x) $$

where $x_t$ is the unique element on $\partial K_t$ that is on the line
through $K_T$ and x. Since $\Bbb P \{(x_1,...,x_n)| x_t \notin
[x_1,...,x_n] \} $ is concentrated  for large n near the boundary
$\partial K$ we get

$$\lim_{n \to \infty} \frac{vol_d(K)-\Bbb
E(K,n)}{(\frac{vol_d(K)}{n})^{\frac{2}{d+1}}}=$$

$$\lim_{n \to \infty}(\frac{n}{vol_d(K)})^{\frac{2}{d+1}}
\int_{\partial K} \int_0^{ \frac{\log n}{n}} \frac{\Bbb P \{(x_1,....,
x_n) \mid x_t \notin [x_1,....,x_n] \} } {vol_{d-1}(P_{N(x_t)}
((-x_{t}+K) \cap (x_{t}-K)))}  \frac {\Vert x_t \Vert^d}{\Vert x \Vert
^d}  \frac {<x,N(x)>}{<x_t,N(x_t)>} dt d\mu $$

Then we apply Lebesgue's convergence theorem and obtain

$$ \int_{\partial K} \lim_{n \to
\infty}(\frac{n}{vol_d(K)})^{\frac{2}{d+1}}\int_0^{\frac{\log n}{n}}
\frac{\Bbb P \{(x_1,...., x_n) \mid x_t \notin [x_1,....,x_n] \} }
{vol_{d-1}(P_{N(x_t)} ((-x_{t}+K) \cap (x_{t}-K)))}  \frac {\Vert x_t
\Vert^d}{\Vert x \Vert ^d}  \frac {<x,N(x)>}{<x_t,N(x_t)>} dt d\mu $$

The hypothesis of Lebesgue's convergence theorem is fulfilled since the
following function dominates the functions under the integral: For
every $x \in \partial K$ let r(x) be the largest radius so that

$$                                B_2^d(x-r(x)N(x),r(x)) \subseteq
K                      \tag 1       $$

r(x) may be 0, e.g. if N(x) is not unique. The functions under the
integral are uniformly smaller than a constant times

$$
r(x)^{-\frac{d-1}{d+1}}
$$

which is integrable on $\partial K$. Then we show that the expression
under the integral converges to $r(x)^{1/(d+1)}$ times an appropriate
constant.  \vskip 1.5cm

\specialhead 3. Proof Of Theorem 1  \endspecialhead

\proclaim{\smc Lemma 2} Let K be a convex body in $R^d$, f a continuous
function on K , and $t \in (0,T)$. Then we have

$$ \frac{d}{dt} \int_{K_t} f(x)dx = -\int_{\partial K_t}
\frac{f(x)}{vol_{d-1}(P_{N(x)} ((-x+K) \cap (x-K)))}
d\mu_{t}(x)                                                        $$

where $\mu_{t}$ is the surface measure on $\partial K_{t}$.
\endproclaim \vskip 1cm For f(x) identical to 1 this is an unpublished
result of Schmuckenschl\"ager [Schm 2]. A similar result for convex
floating bodies instead of convolution bodies can be found in
[Sch\"u-$W_2$].  \vskip 1cm

\proclaim{\smc Lemma 3} Let K be a convex body in $R^d$ and $K_t$, $t
\in [0,T]$, the convolution bodies. Then we have for all $t_0 \in
[0,T]$

$$  vol_d(K)- \Bbb
E(K,n)=
$$

$$ \int_{\partial K} \int_0^{t_0} \frac{\Bbb P \{(x_1,...., x_n) \mid
x_t \notin [x_1,....,x_n] \} } {vol_{d-1}(P_{N(x_{t})} ((-x_{t}+K) \cap
(x_{t}-K)))}  \frac {\Vert x_t \Vert^d}{\Vert x \Vert ^d}  \frac
{<x,N(x)>}{<x_t,N(x_t)>} dt d\mu(x) $$

$$+ \int_{t_0}^T \int_{\partial K_t} \frac {\Bbb P\{(x_1,...,x_n) \mid
x \notin [x_1,...,x_n]\}}{vol_{d-1}(P_{N(x)} ((-x+K) \cap (x-K)))}
d\mu_t(x)dt           $$

where $\mu$ and $\mu_{t}$ are the surface measures on $\partial K$ and
$\partial K_{t}$ respectively and $\{x_t\}=\partial K \cap [K_T,x]$.
\endproclaim \vskip 1cm

\demo{\smc Proof}

$$ vol_d(K)- \Bbb E(K,n) = \int_{K \times \hdots \times K} vol_d(K
\setminus [x_1,...,x_n]) d \Bbb P    $$

$$ = \int_{K \times \hdots \times K} \int_{K} \chi_{K \setminus
[x_1,...,x_n]} dx d\Bbb P                $$

$$ = \int_{K} \int_{K \times \hdots \times K} \chi_{K \setminus
[x_1,...,x_n]} d\Bbb P dx               $$

$$ = \int_{K} \Bbb P \{(x_1,...,x_n)| x \notin [x_1,...,x_n]\}
dx                                        $$

Since

$$  \int_{K_t} \Bbb P \{(x_1,...,x_n)| x \notin [x_1,...,x_n] \}
dx                                      $$

is a bounded, continuous, decreasing function on $[0,T]$ it is
absolutely continuous. We get

$$ vol_d(K)- \Bbb E(K,n)= -\int_0^T \frac{d}{dt} \int_{K_t} \Bbb P
\{(x_1,...,x_n)| x \notin [x_1,...,x_n] \}
dxdt
$$

Since $\Bbb P \{(x_1,...,x_n)| x \notin [x_1,...,x_n] \}$ is a
continuous function of x we get by Lemma 2

$$ vol_d(K)- \Bbb E(K,n)= \int_0^T \int_{\partial K_t} \frac{\Bbb P
\{(x_1,...,x_n)| x \notin [x_1,...,x_n] \} }{vol_{d-1}(P_{N(x)}((-x+K)
\cap (x-K)))} d\mu_{t}(x)
dt                                               $$

\qed \enddemo \vskip 1cm

\proclaim{\smc Lemma 4} Let cap(r,$\Delta$) be a cap of height $\Delta$
of a d-dimensional Euclidean sphere with radius r.  Then we have

$$  2(2- \frac{\Delta}{r})^{\frac{d-1}{2}}
\frac{vol_{d-1}(B_{2}^{d-1})}{d+1} \Delta^{\frac{d+1}{2}}
r^{\frac{d-1}{2}}
$$

$$  \leq vol_d(cap(r,\Delta))
\leq
$$

$$ 2^{\frac{d+1}{2}} \frac{vol_{d-1}(B_{2}^{d-1})}{d+1}
\Delta^{\frac{d+1}{2}} r^{\frac{d-1}{2}}           $$ \endproclaim
\vskip 1cm

\proclaim{\smc Lemma 5} Let K be a convex body in $\Bbb R^d$. Then
there are constants $c,c' > 0$ so that we have for all $x \in \partial
K$ and for all $r > 0$ with $B_{2}^{d}(x-rN(x),r) \subseteq K$

$$ vol_{d-1}(P_{N(x_t)}((-x_{t}+K) \cap (x_{t}-K))) \geq
   \cases c(tr)^{\frac{d-1}{d+1}} &\text{if $0 \leq t \leq c'r^{d}$}
   \\ ct^{\frac{d-1}{d}} &\text{if $c'r^{d} \leq t \leq
   T$}                   \endcases                     $$

where $\{x_{t}\}= \partial K_{t} \cap [x,K_{T}]$.  \endproclaim \vskip
1cm

\demo{\smc Proof} We may assume that $K_T$ coincides with the origin.
By convexity and the fact that $K_T$ is an interior point we get that
there is $c_1 > 0$ so that we have for all $x \in \partial K$

$$ <\frac{x}{\parallel x \parallel}, N(x)> \geq
c_1                                     \tag 2   $$

Now we choose

$$ c= \frac{vol_{d-1}(B_2^{d-1})}{(d+1)vol_d(B_2^d)} min\{(1- \root
\of{1-c_1^2})^{\frac{d+1}{2}},
   (1- \root \of{\frac{3}{4}})^{\frac{d+1}{2}}
   \}                                               \tag 3    $$

Then we have for all $t \in [0,cr^d vol_d(B_2^d)]$

$$  \parallel x-x_t \parallel  \leq r<\frac{x}{ \parallel x \parallel},
N(x)>             \tag 4       $$

$$  \frac{1}{2}  \leq
|<N(x),N(x_t)>|                                                 \tag
5         $$

Geometrically (4) means the following: Let z be the midpoint of the
interval $[K_T,x] \cap B_2^d(x-rN(x),r)$. Then $x_t \in [z,x]$ (see
figure 1).

\midinsert \vspace{12cm} \botcaption{Figure 1} \endcaption \endinsert

We verify (4). Since $B_2^d(x-rN(x),r) \subseteq K$ we have

$$ vol_d((-z+K) \cap (z-K)) \geq vol_d((-z+B_2^d(x-rN(x),r)) \cap
(z-B_2^d(x-rN(x),r)))         $$

The last expression equals twice the volume of a cap whose height is
greater than

$$  \Delta=r(1- \root \of{1-c_1^2})            $$

By Lemma 4 and (3) we get

$$ vol_d((-z+K) \cap (z-K)) \geq 4(2- \frac{\Delta}{r})^{\frac{d-1}{2}}
\frac{vol_{d-1}(B_2^{d-1})}{d+1}
  \Delta^{\frac{d+1}{2}}
  r^{\frac{d-1}{2}}                                               $$

$$ \geq 4cr^d vol_d(B_2^d)  \geq 4t >t              $$

Thus z is an interior point of $K_t$ and $x_t \in [z,x]$. Moreover, $
\parallel z-x \parallel = r<\frac{x}{\parallel x \parallel}, N(x)> $.
\par Now we check (5). We get by (4) and figure 1 that $x_t$ has to be
in the shaded area of figure 2.

\midinsert \vspace{9cm} \botcaption{Figure 2} \endcaption \endinsert

Assume that (5) is not true. Then it follows that the radius of the
sphere $H(x_t,N(x_t)) \cap B_2^d(x-rN(x),r)$ is greater than r and that
$H(x_t,N(x_t))$ cuts off a cap of height greater than $(r(1- \root \of
{\frac{3}{4}}))$. By Lemma 4 we get that the volume of this cap is
greater than

$$ 2(2- \frac{\Delta}{r} )^{\frac{d-1}{2}}
\frac{vol_{d-1}(B_2^{d-1})}{d+1} r^d(1- \root
   \of{\frac{3}{4}})^{\frac{d+1}{2}}
   $$

Therefore we have for the center w of the sphere $B_2^d(x-rN(x),r) \cap
H(x_t,N(x_t))$ that

$$ vol_d((-w+K) \cap (w-K)) \geq 2 vol_d(B_2^d(x-rN(x),r) \cap
H^{-}(x_t,N(x_t)))             $$

$$ \geq 4 \frac{vol_{d-1}(B_2^{d-1})}{d+1} (1- \root
\of{\frac{3}{4}})^{\frac{d+1}{2}} r^d
   \geq 4t
   >t
   $$

This means that w is an interior point of $K_t$. On the other hand, z
is an element of the supporting hyperplane $H(x_t,N(x_t))$ to $K_t$.
This gives a contradiction and we conclude that (5) is valid. \par We
denote $\Theta = \arccos(<\frac{x}{\parallel x \parallel}, N(x)>)$.

\midinsert \vspace{13cm} \botcaption{Figure 3} \endcaption \endinsert

>From figure 3 it follows that the distance of $x_t$ to the boundary of
$B_2^d(x-rN(x),r)$ equals

$$ (\cos(\Theta)-\sin(\Theta) \cot(\frac{\pi}{2} - \frac{\alpha}{2}))
\parallel x-x_t \parallel    $$

By figure 1 we have $\alpha \leq \frac{\pi}{2} - \Theta$ so that the
above expression is larger than

$$ (\cos(\Theta)- \sin(\Theta) \cot(\frac{\pi}{4}+ \frac{\Theta}{2}))
\parallel x-x_t \parallel     \tag 6      $$

Please note that by (2) there is $\epsilon > 0$ so that we have for all
$x \in \partial K$

$$   \cos(\Theta) - \sin(\Theta) \cot(\frac{\pi}{4}+\frac{\Theta}{2})
\geq \epsilon             \tag 7     $$

We assume now that

$$  \parallel x-x_t \parallel \geq
\frac{1}{\epsilon}\frac{t^{\frac{2}{d+1}}}{r^{\frac{d-1}{d+1}}}
    (\frac{d+1}{vol_{d-1}(B_2^{d-1})})^{\frac{2}{d+1}}
    \tag 8     $$

Then

$$ (-x_t+K) \cap (x_t-K) \supseteq (-x_t+B_2^d(x-rN(x),r)) \cap
(x_t-B_2^d(x-rN(x),r))       $$

has a volume greater than twice the volume of a cap of a Euclidean ball
of radius r and height (6). By Lemma 4, (6), and (7) we get as above

$$ vol_d((-x_t+K) \cap (x_t-K)) \geq 4t           $$

which cannot be true. \par Therefore (8) does not hold. We deduce that
the distance between the two parallel hyperplanes $H(x,N(x))$ and
$H(2x_t-x,-N(x))$ is less than twice the right hand expression of (8).
Moreover,

$$ (-x_t+K) \cap (x_t-K) \subseteq H^{+}(x,N(x)) \cap
H^{+}(2x_t-x,-N(x))      $$

where both half spaces are chosen so that $x_t$ is contained in them.
Therefore there must be a hyperplane H parallel to $H(x,N(x))$ so that

$$ vol_{d-1}((-x_t+K) \cap (x_t-K) \cap
H)                                        $$

$$   \geq \frac{\epsilon}{2} (rt)^{\frac{d-1}{d+1}}
   (\frac{vol_{d-1}(B_2^{d-1})}{d+1})^{\frac{2}{d+1}}
   $$

By (5) we get that the same inequality holds for

$$ P_{N(x_t)}((-x_t+K) \cap (x_t-K))                $$

with another constant.

\par Now we consider the case

$$ cr^d vol_d(B_2^d) \leq t \leq T       $$

Since K is compact there are $r_1, r_2 > 0$ so that

$$ B_2^d(K_T,r_1) \subseteq K \subseteq
B_2^d(K_T,r_2)                         \tag 9   $$

Suppose that

$$ vol_{d-1}(P_{N(x_t)}((-x_t+K) \cap (x_t-K)))) \leq
ct^{\frac{d-1}{d}}        \tag 10  $$

with $ c< \frac{r_1}{4r_2} vol_d(B_2^d)^{\frac{1}{d}}$. This implies
that

$$ vol_{d-1}(((-x_t+K) \cap (x_t-K)) \cap H(x_t,N(x_t)) \leq
ct^{\frac{d-1}{d}}         $$

and since $(-x_t+K) \cap (x_t-K)$ is symmetric with respect to $x_t$
all other $d-1$ dimensional sections of $(-x_t+K) \cap (x_t-K)$ that
are parallel to $H(x_t,N(x_t))$ have a smaller $d-1$ dimensional
volume. Therefore there must be a non-empty section whose distance to
$H(x_t,N(x_t))$ is at least $\frac{1}{2c} t^{\frac{1}{d}}$. This means
that there is $z \in K$ with (see figure 4)

$$ d(z,H(x_t,N(x_t)) \geq \frac{1}{2c}
t^{\frac{1}{d}}                        \tag 11   $$

Let y be the unique point in the intersection

$$ [K_T,z] \cap H(x_t,N(x_t))        $$

We show that y is an interior point of $K_t$ which contradicts the fact
that $y \in H(x_t,N(x_t))$.  The sphere

$$ [z, B_2^d(K_T,r_1)] \cap H(y, \frac{z-y}{\parallel z-y
\parallel})       $$

has a radius larger than

$$ \frac{1}{2c} \frac{r_1}{r_2} t^{\frac{1}{d}}           $$

\midinsert \vspace{15cm} \botcaption{Figure 4} \endcaption \endinsert

This follows from (11). Thus we find that

$$ [z, B_2^d(K_T, r_1)] \supseteq B_2^d(y, \frac{1}{4c} \frac{r_1}{r_2}
t^{\frac{1}{d}})   $$

and as above that y is an interior point of $K_t$.

\qed \enddemo \vskip 1cm

\proclaim{\smc Lemma 6} [Sch\"u-W1] Let K be a convex body in $\Bbb
R^d, \alpha \in (0,1)$ and r(x) as defined by (1). Then we have

$$   \int_{\partial K} r(x)^{-\alpha} d\mu(x) < \infty
$$

where $\mu$ is the surface measure on $\partial K$.  \endproclaim
\vskip 1cm

\proclaim{\smc Lemma 7} [B-L] Let K be a convex body in $\Bbb R^d$ and
let $x \in \partial K_t$. Then we have

$$ \Bbb P \{(x_1,...,x_n)| x_t \notin [x_1,...,x_n] \}
 \leq 2 \sum_{i=0}^{d-1} \binom{n}{i} (\frac{t}{2vol_d(K)})^{i} (1-
 \frac{t}{2vol_d(K)})^{n-i}   $$ \endproclaim \vskip 1cm

\proclaim{\smc Lemma 8} Let K be a convex body in $\Bbb R^d$, $t_0 \in
(0,T]$ and $\mu_t$ the surface measure on $\partial K_t$. Then we have

$$ \lim_{n \to \infty} n^{\frac{2}{d+1}} \int_{t_0}^{T} \int_{\partial
K_t} \frac{\Bbb P \{ (x_1,...,x_n)| x \notin [x_1,...,x_n] \} }
     {vol_{d-1}(P_{N(x)}((-x+K) \cap (x-K)))} d\mu_{t}(x) dt
     =0               $$ \endproclaim \vskip 1cm

\demo{\smc Proof} Since $t \in [t_0,T]$ and $t_0 > 0$ there is a
constant $c>0$ so that we have for all $x \in \partial K_t$

$$   vol_{d-1}(P_{N(x)}((-x+K) \cap (x-K)) \geq c    $$

By Lemma 7 we get

$$ n^{\frac{2}{d+1}} \int_{t_0}^{T} \int_{\partial K_t} \frac{\Bbb P \{
(x_1,...,x_n)| x \notin [x_1,...,x_n] \} }
     {vol_{d-1}(P_{N(x)}((-x+K) \cap (x-K)))} d\mu_{t}(x)
     dt                $$

$$ \leq \frac{2}{c} n^{\frac{2}{d+1}} \sum_{i=0}^{d-1} \binom{n}{i}
\int_{t_0}^{T} \int_{\partial K_t} (\frac{t}{2vol_d(K)})^{i} (1-
\frac{t}{2vol_d(K)})^{n-i}d\mu_{t}(x) dt   $$

Since $vol_{d-1}(\partial K_t) \leq vol_{d-1}(\partial K)$ we get that
the last expression is smaller than

$$  \frac{2}{c} n^{\frac{2}{d+1}} vol_{d-1}(\partial K)
 \sum_{i=0}^{d-1} \binom{n}{i} \int_{t_0}^{T} (\frac{t}{2vol_d(K)})^{i}
 (1- \frac{t}{2vol_d(K)})^{n-i}dt   $$

$$ \leq \frac{4}{c} n^{\frac{2}{d+1}} vol_{d-1}(\partial K) vol_{d}(K)
\sum_{i=0}^{d-1} \binom{n}{i} \frac{ \Gamma(i+1)
\Gamma(n-i+1)}{\Gamma(n+2)}  $$

$$= \frac{4}{c} vol_{d-1}(\partial K) vol_{d}(K) n^{\frac{2}{d+1}}
\frac{d}{n+1}  $$

$$ \leq \frac{4}{c} vol_{d-1}(\partial K) vol_{d}(K) d
n^{-\frac{d-1}{d+1}}      $$

\qed \enddemo \vskip 1cm

\proclaim{\smc Lemma 9} Let K be a convex body in $\Bbb R^d$, let $K_T$
be the origin, and $0 \leq t_1 \leq t_2 < T$.  Then there is a constant
$c>0$ so that we have for all $n \in \Bbb N$

$$
 \int_{t_1}^{t_2} \frac{\Bbb P \{(x_1,...., x_n) \mid x_t \notin
[x_1,....,x_n] \} } {vol_{d-1}(P_{N(x_t)} ((-x_{t}+K) \cap
(x_{t}-K)))}  \frac {\Vert x_t \Vert^d}{\Vert x \Vert ^d}  \frac
{<x,N(x)>}{<x_t,N(x_t)>} dt   $$

\leftline{(12)}

$$ \leq c r(x)^{-\frac{d-1}{d+1}} \int_{\frac{t_1}{2vol_d(K)}}^{\frac{
t_2}{2vol_d(K)}} \sum_{i=0}^{d-1} \binom{n}{i}  s^{i- \frac{d-1}{d+1}}
(1-s)^{n-i}
ds                                                             $$

where r(x) is defined by (1).  \endproclaim \vskip 1cm

\demo{\smc Proof} We have $\parallel x_t \parallel \leq \parallel x
\parallel$ and we have a constant $c>0$ so that we have for all $x \in
\partial K$ $<x_t,N(x_t)> \geq c$ because $t_2 < T$. Thus it is enough
to estimate

$$
 \int_{t_1}^{t_2} \frac{\Bbb P \{(x_1,...., x_n) \mid x_t \notin
[x_1,....,x_n] \} } {vol_{d-1}(P_{N(x_t)} ((-x_{t}+K) \cap (x_{t}-K)))}
 dt   $$

By Lemma 5 this is smaller than

$$ \int_{t_1}^{c'r^d} c(r(x)t)^{-\frac{d-1}{d+1}} \Bbb P
\{(x_1,...,x_n)| x_t \notin [x_1,...,x_n]
\}dt                                             $$

$$ +\int_{c'r^d}^{t_2} ct^{-\frac{d-1}{d}} \Bbb P \{(x_1,...,x_n)| x_t
\notin [x_1,...,x_n] \}dt
$$

Since we have for $t \in [c'r^d,t_2]$ that

$$ t^{-\frac{d-1}{d}} = t^{-\frac{d-1}{d(d+1)}} t^{-\frac{d-1}{d+1}}
  \leq (c'r(x))^{-\frac{d-1}{d(d+1)}} t^{-\frac{d-1}{d+1}}   $$

we can estimate the above expression by

$$ cr(x)^{-\frac{d-1}{d+1}} \int_{t_1}^{t_2} t^{-\frac{d-1}{d+1}} \Bbb
P \{(x_1,...,x_n)| x_t \notin [x_1,...,x_n]
\}dt                                             $$

where c is a new constant. Now it is left to apply Lemma 7.  \qed
\enddemo \vskip 1cm

\proclaim{\smc Lemma 10} Let K be a convex body in $\Bbb R^d$, let
$K_T$ be the origin, and let $t_1 < T$.  Then there is a constant c so
that we have for all $x \in \partial K$ and all $n \in \Bbb N$

$$  \int_0^{t_1} \frac{\Bbb P \{(x_1,...., x_n) \mid x_t \notin
[x_1,....,x_n] \} } {vol_{d-1}(P_{N(x_{t})} ((-x_{t}+K) \cap
(x_{t}-K)))}  \frac {\Vert x_t \Vert^d}{\Vert x \Vert ^d}  \frac
{<x,N(x)>}{<x_t,N(x_t)>} dt \leq c r(x)^{-\frac{d-1}{d+1}} $$

\endproclaim \vskip 1cm

\demo{\smc Proof} By Lemma 9 we get

$$  \int_0^{t_1} \frac{\Bbb P \{(x_1,...., x_n) \mid x_t \notin
[x_1,....,x_n] \} } {vol_{d-1}(P_{N(x_t)} ((-x_{t}+K) \cap
(x_{t}-K)))}  \frac {\Vert x_t \Vert^d}{\Vert x \Vert ^d}  \frac
{<x,N(x)>}{<x_t,N(x_t)>} dt     $$

$$ \leq cr(x)^{-\frac{d-1}{d+1}} n^{\frac{2}{d+1}}  \sum_{i=0}^{d-1}
\binom{n}{i} \frac{\Gamma(i+1-\frac{d-1}{d+1})
\Gamma(n+1-i)}{\Gamma(n+2-\frac{d-1}{d+1})}  $$

Since

$$ \lim_{k \to \infty} \frac{\Gamma(k+\frac{2}{d+1})}{\Gamma(k)}
k^{-\frac{2}{d+1}}=1   $$

we can estimate the last expression by

$$  cr(x)^{-\frac{d-1}{d+1}} $$

where c is a new constant that does not depend on n and x.  \qed
\enddemo \vskip 1cm

\proclaim{\smc Lemma 11} We have

$$ \lim_{n \to \infty} n^{\frac{2}{d+1}} \int_{\frac{\log n}{n}}^{1}
\sum_{i=0}^{d-1} \binom{n}{i} s^{i-\frac{d-1}{d+1}}(1-s)^{n-i} ds=0
$$ \endproclaim \vskip 1cm

\proclaim{\smc Lemma 12} Let K be a convex body in $\Bbb R^d$, let
$K_T$ be the origin, and let $t_1 < T$.  Then we have for all $x \in
\partial K$ with $r(x)>0$

$$ \lim_{n \to \infty} n^{\frac{2}{d+1}}\int_{\frac{\log n}{n}}^{t_1}
\frac{\Bbb P \{(x_1,...., x_n) \mid x_t \notin [x_1,....,x_n] \} }
{vol_{d-1}(P_{N(x_t)} ((-x_{t}+K) \cap (x_{t}-K)))}  \frac {\Vert x_t
\Vert^d}{\Vert x \Vert ^d}  \frac {<x,N(x)>}{<x_t,N(x_t)>} dt=0 $$
\endproclaim \vskip 1cm

Proof. The result follows from Lemmata 9 and 11.  \qed \vskip 1cm

\proclaim{\smc Lemma 13}[Wie]

$$ \lim_{n \to \infty} n^{\frac{2}{d+1}}( vol_{d}(B_2^{d}(0,r))- \Bbb
E(B_2^{d}(0,r),n))=       $$

$$ \frac{(d^2+d+2)(d^2+1)}{2(d+3)(d+1)!}
   ((d+1) \frac{vol_d(B_2^d)}{vol_{d-1}(B_2^{d-1})})^{\frac{2}{d+1}}
   \Gamma(\frac{d^2+1}{d+1}) vol_{d-1}(\partial B_2^d(0,r))      $$

\endproclaim \vskip 1cm

By an affine transform we can change the indicatrix of Dupin into a
Euclidean sphere or a cylinder with a sphere as its base.  \vskip 1cm

\proclaim{\smc Lemma 14} Let K be a convex body in $\Bbb R^d$ with $0
\in \partial K$ and $N(0)=(0,...,0,-1)$. Suppose that $\partial K$ is
twice differentiable in the generalized sense at 0. \par (i) If the
indicatrix of Dupin at 0 is a d-2 dimensional sphere with radius $\root
\of \rho$, then there is a $t_0 > 0$ and a monotone , increasing
function $\psi$ on $\Bbb R^+$ with $\lim_{t \to 0} \psi(t)=\psi(0)=1 $
so that we have for all $t \in (0,t_0]$

$$ \{(\frac{x_1}{\psi(t)},...,\frac{x_{d-1}}{\psi(t)},t)|
    x \in B_2^d((0,...,0,\rho),\rho) \text{\ and\ } x_n=t
    \}             $$

$$ \subseteq K \cap H(-tN(0),N(0))                             $$

$$ \subseteq \{(\psi(t)x_1,...,\psi(t)x_{d-1},t)|
    x \in B_2^d((0,...,0,\rho),\rho) \text{\ and\ } x_n=t
    \}             $$

(ii) If the indicatrix of Dupin at 0 is a d-2 dimensional cylinder with
radius $\root \of \rho$, i.e.

$$  \Bbb R^{k-1} \times \partial B_2^{d-k}(0,\root \of \rho)      $$

then there is a function $\Phi$ on $\Bbb R^+$ so that for every
$\epsilon >0$ there is a $t_0 > 0$ so that $\lim_{t \to 0} \frac{\root
\of t}{\Phi(t)}=0$ and $\frac{\root \of t}{\Phi(t)}$ is increasing on
$\Bbb R^+$ and so that we have for all $t \in (0,t_0]$

$$ \{(y,x,t)| (x,t) \in B_2^{d-k+1}((0,...,0,\rho-\epsilon), \rho -
\epsilon)\text{\ and\ }
   y \in [-\Phi(t),\Phi(t)]^{k-1}   \}
   $$

$$ \subseteq K \cap H(-tN(0),N(0))                                 $$

\endproclaim \vskip 1cm

\proclaim{\smc Lemma 15} Let K be a convex body in $\Bbb R^d$, $c>2$,
 and $x \in \partial K$ such that $\kappa (x) >0$.  Then there is
$t_{c} >0$ so that we have for all $t \in (0,t_{c}]$: We have for the
hyperplane H whose normal coincides with N(x) and that satisfies
$vol_{d}(K \cap H^{-})=c^{d+1}t$ and for all $n \in \Bbb N$ with $n
\geq \frac{vol_{d}(K)} {ct}$

$$ |\Bbb P\{(x_1,...,x_n)| x_t \notin [\{x_1,...,x_n\} \cap H^{-}]\}
   - \Bbb P\{(x_1,...,x_n)| x_t \notin
   [x_1,...,x_n]\}|                    $$

$$ < 2^{d-1} e^{-c_1 c}           $$

where $c_1$ is a constant that depends on d only.

\endproclaim \vskip 1cm

\demo{\smc Proof} We show that

$$ \Bbb P\{(x_1,...,x_n)| x_t \notin [\{x_1,...,x_n\} \cap H^{-}]\}
   \text{\ and\ } x_t \in [x_1,...,x_n]\} \leq 2^{d-1} e^{-c_1
   c}                             $$

If we have

$$ x_t \notin [\{x_1,...,x_n\} \cap H^{-}]\}\text{\ and\ } x_t \in
[x_1,...,x_n] $$

then there is $y \in H^{+} \cap K$ so that

$$ [y,x_t] \cap [\{x_1,...,x_n\} \cap H^{-}] = \emptyset
\tag 13 $$

For the following argument let us assume that $N(x)=e_1$ and $x_t=0$.
Moreover, since $\kappa (x)>0$ we may assume that the indicatrix of
Dupin at x is a Euclidean sphere and by Lemma 14 $\partial K$ can be
approximated arbitrary well at x by a sphere if we choose the height of
the cap or correspondingly $t_{c}$ sufficiently small. We assume for
the following arguments that $K \cap H^{-}$ is a cap of a sphere. Later
we shall see that we control the error by choosing $t_{c}$ sufficiently
small. We consider the following sets (figure 5)

$$ corn_{\Theta}=K \cap H^{-} \cap H^{+}(x_t,N(x)) \cap
   \{\bigcap_{i=2}^{d}H^{-}(x_t,e_1+\Theta_{i}\lambda e_i)\}      $$

where $\Theta_2,...,\Theta_d= \pm 1$. We have $2^{d-1}$ sets and they
are
 best described as corner sets. $\lambda$ is chosen so that

$$ H^{-}(x_t, e_1+\Theta_{i} \lambda e_i) \cap H^{+}  $$

consists of exactly one point.

\midinsert \vspace{8cm} \botcaption{Figure 5} \endcaption \endinsert

By the height $h_1$ of $corn_{\Theta}$ we understand the minimal
distance of $H(x_t,N(x))$ and a parallel hyperplane so that
$corn_{\Theta}$ lies between them. We get

$$ vol_{d}(corn_{\Theta}) \geq 2^{-d+1} \frac{h_1}{d}
   vol_{d-1}(K \cap H(x_t,N(x))                               $$

Let $h_2$ denote the height of the cap $K \cap H^{-}(x_t,N(x))$ and
$h_3$ the height of $K \cap H^{-}$. By Lemma 4 we get that

$$ 2c^{d+1}=\frac{vol_{d}(K \cap H^{-})}{vol_{d}(K \cap
H^{-}(x_{t},N(x)))}
  \leq 2^{\frac{d+1}{2}}
  (\frac{h_3}{h_2})^{\frac{d+1}{2}}                  $$

or

$$ h_{2} c^{2} \leq 2 h_{3}  $$

The height $h_1$ is of the order $\frac{h_3}{d}$. Altogether we get

$$ vol_{d}(corn_{\Theta}) \geq 2^{-d+1} c_0 c^2 h_2
   vol_{d-1}(K \cap H(x_t,N(x)) \geq c_2 c^{2}
   t                                    \tag 14  $$

where $c_2$ is a constant depending only on d.  We have by (13)

$$ \{(x_1,...,x_n)| x_t \notin [\{x_1,...,x_n\} \cap H^{-}]
\text{\ and\ } x_t
   \in [x_1,...,x_n]\}
   \subseteq                                                  $$

$$ \{(x_1,...,x_n)| \exists H_{x_t} : x_t \in H_{x_t},
   H^{-}_{x_t} \cap K \cap H^{+} \ne \emptyset \text{\ and\ }
   H^{+}_{x_t} \supset \{x_1,...,x_n\} \cap H^{-}
   \}                                          $$

Indeed, by the theorem of Hahn-Banach there is a hyperplane $H_{x_t}$
separating the convex sets $[\{x_1,...,x_n\} \cap H^{-}]$ and the ray
$\{x_{t}+ \lambda (y-x_{t})| \lambda \in \Bbb R \}$. By (13) they are
disjoint. We may assume that at least one point of the ray is an
element of $ H_{x_t}$. So $x_t$ is also an element of $H_{x_t}$. Let
$H^{-}_{x_t}$ be the halfspace containing y, then $H^{-}_{x_t} \cap K
\cap H^{+}$ contains y and is not empty. \par Such a halfspace always
contains one of the corner sets $corn_{\Theta}$.  This follows since we
have in $\Bbb R^d$ for a hyperplane passing through the origin: The
corresponding halfspaces contain at least one $2^d$-tant. Therefore we
get

$$  \{(x_1,...,x_n)| \exists H_{x_t} :
   H^{-}_{x_t} \cap K \cap H^{+} \ne \emptyset \text{\ and\ }
   H^{+}_{x_t} \supset \{x_1,...,x_n\} \cap H^{-}
   \}                                          $$

$$\subseteq \bigcup_{\Theta} \{(x_1,...,x_n)|\{x_1,...,x_n\}
   \subseteq K \setminus corn_{\Theta} \}
   $$

And consequently we get by (14)

$$\Bbb P \{(x_1,...,x_n)| x_t \notin [\{x_1,...,x_n\} \cap H^{-}]
\text{\ and\ } x_t
   \in
   [x_1,...,x_n]\}
   $$

$$ \leq 2^{d-1} \Bbb P\{(x_1,...,x_n)|\{x_1,...,x_n\} \subset K
\setminus corn_{\Theta} \}   $$

$$ \leq 2^{d-1}(1- \frac{c_2 c^2 t}{vol_d(K)})^{n}
   \leq 2^{d-1} \exp(-n \frac{c_2 c^2
   t}{vol_d(K)})                            $$

By the assumption on n we get that the last expression is smaller than

$$ 2^{d-1} \exp(-c_2c)  $$

This argument also works if the volumes of the considered sets differ
by a small error. Therefore the proof also goes through if K is not a
sphere but can be approximated arbitrarily well by a sphere at the
point x.  \qed \enddemo \vskip 1cm

\proclaim{\smc Lemma 17} Let K be a convex body in $\Bbb R^d$ and B a
Euclidean ball of the same volume. Let $x \in \partial K$ and $z \in
\partial B$ and assume that $\kappa(x) >0$. Then for every $\epsilon >
0$ there is $t_{\epsilon} > 0$ so that we have for all $t \in
(0,t_{\epsilon}]$ and all $n \in \Bbb N$ with $n \geq 2d$

$$| \Bbb P_{K} \{(x_1,...,x_n)| x_t \in [x_1,...,x_n] \}
     - \Bbb P_{B} \{(z_1,...,z_n)| z_t \in [z_1,...,z_n] \}| <
     \epsilon       \tag 15   $$

\endproclaim \vskip 1cm

\demo{\smc Proof} We show first that there is $c>1$ so that we have
(15) whenever $n \leq \frac{vol_d(K)} {ct}$ or $n \geq c
\frac{vol_d(K)}{t}$. As c we can choose a number satisfying

$$ c\geq max\{\frac{1}{\epsilon},d \}     $$

$$ (d4c^{d+2})^{d+1} e^{-\frac{c}{2}} < \epsilon            \tag 16  $$

$$ 2^{d-1} e^{-c_1 c}< \epsilon     $$

where $c_1$ is the constant introduced in Lemma 15.  \par We consider
the case $n \leq \frac{vol_d(K)}{ct}$. Since the curvature at x is
strictly positive we may assume that the indicatrix of Dupin is a
sphere. If we choose $t_{\epsilon}$ small enough then by Lemma 14 there
is a hyperplane H through $x_t$,$ 0 < t \leq t_{\epsilon}$ , so that
$vol_d(K \cap H^{-}) \leq t$. Therefore we get

$$ 1 \geq \Bbb P_{K} \{(x_1,...,x_n)| x_t \notin [x_1,...,x_n] \}
 \geq \Bbb P_{K} \{(x_1,...,x_n)| \{x_1,...,x_n\} \subseteq K \cap
 H^{+} \}  $$

$$   \geq  (1- \frac{t}{vol_d(K)})^{n} \geq (1- \frac{1}{cn})^{n} \geq
1- \frac{1}{c}   $$

The same estimate holds for $\Bbb P_{B}$ and we get (15). Now we
consider the case $n \geq 2c \frac{vol_d(K)}{t}$. By Lemma 7 we get

$$0 \leq \Bbb P \{(x_1,...,x_n)| x_t \notin [x_1,...,x_n] \}    $$

$$ \leq 2 \sum_{i=0}^{d-1} \binom{n}{i} (\frac{t}{2vol_d(K)})^{i} (1-
\frac{t}{2vol_d(K)})^{n-i}   $$

The function $s^{i}(1-s)^{n-i}$ attains its maximum at $\frac{i}{n}$.
Since $i < d$, $d \leq c$, and $2d \leq n$ we get that the latter
expression is less than

$$ \sum_{i=0}^{d-1} \binom{n}{i} (\frac{c}{n})^{i} (1-
\frac{c}{n})^{n-i} \leq
    d c^d exp(- \frac{c}{2}) \leq
    \epsilon                                  $$

The same holds for $\Bbb P_{B}$ and we get (15) again.   \par Now we
consider the case $ \frac{vol_{d}(K)}{ct} \leq n \leq c
\frac{vol_{d}(K)}{t}$.  By triangle inequality we get

$$ | \Bbb P_{K} \{(x_1,...,x_n)| x_t \notin [x_1,...,x_n] \} -
     \Bbb P_{B} \{(z_1,...,z_n)| z_t \notin [z_1,...,z_n] \} | \leq  $$

$$ | \Bbb P_{K} \{(x_1,...,x_n)| x_t \notin [x_1,...,x_n] \}-
     \Bbb P_{K} \{(x_1,...,x_n)| x_t \notin [\{x_1,...,x_n\} \cap
     H^{-}] \}|    $$

$$ +| \Bbb P_{K} \{(x_1,...,x_n)| x_t \notin [\{x_1,...,x_n\} \cap
H^{-}] \} -
     \Bbb P_{B} \{(z_1,...,z_n)| z_t \notin [\{z_1,...,z_n\} \cap
     \tilde H^{-} ] \} |  $$

$$ +| \Bbb P_{B} \{(z_1,...,z_n)| z_t \notin [z_1,...,z_n] \} -
      \Bbb P_{B} \{(z_1,...,z_n)| z_t \notin [\{z_1,...,z_n\} \cap
      \tilde H^{-} ] \} |  $$

where H and $\tilde H$ are hyperplanes whose normals coincide with N(x)
and N(z) respectively and $vol_{d}(K \cap H^{-})=vol_d(B \cap \tilde
H^{-}) =c^{d+1}t$.  The first and third summands of the latter
expression can be estimated by Lemma 15.  We estimate now the second
summand. Again, we may assume that the indicatrix of Dupin at $x \in
\partial K$ is a Euclidean sphere. Moreover, we may assume that the
radius of the indicatrix equals the radius of B. This is done by a
volume preserving, affine transform. We have

$$ \Bbb P_{K} \{(x_1,...,x_n)| x_t \notin [\{x_1,...,x_n\} \cap H^{-}]
\}=   $$

\leftline{(17)}

$$  \sum_{k=o}^{n} \binom{n}{k} (\frac{vol_{d}(K \cap
H^{-})}{vol_{d}(K)})^{k}
     (1- \frac{vol_{d}(K \cap H^{-})}{vol_{d}(K)})^{n-k} \Bbb P_{K \cap
     H^{-}} \{(x_1,...,x_k)| x_t \notin [x_1,...,x_k] \}          $$

and the same for $\Bbb P_{B}$. We get by Lemma 7 for $k \geq 4dc^{d+2}
$

$$ \Bbb P_{K \cap H^{-}} \{(x_1,...,x_k)| x_t \notin [x_1,...,x_k]
\}    $$

$$ \leq 2 \sum_{i=0}^{d-1} \binom{k}{i} (\frac{1}{2c^{d+1}})^{i}
    (1- \frac{1}{2 c^{d+1}})^{k-i} \leq 2 d k^{d} e^{-\frac{1}{4} k
 c^{-d-1}}                     $$

The function $s^{d} e^{-as}$ attains its maximum at $\frac{d}{a}$.
Therefore, and because of $4dc^{d+2} \leq k$ the last expression is
smaller than

$$ 2d(d4c^{d+2})^{d} e^{-dc} < \epsilon                              $$

We get the same for $\Bbb P_B$. Therefore we have

$$ | \Bbb P_{K} \{(x_1,...,x_n)| x_t \notin [\{x_1,...,x_n\} \cap
H^{-}] \} - \Bbb P_{B} \{(z_1,...,z_n)| z_t \notin [\{z_1,...,z_n\}
\cap \tilde H^{-} ] \} |  $$

$$ \leq  \sum_{0 \leq k \leq 4dc^{d+2}} \binom{n}{k}
    (\frac{vol_{d}(K \cap H^{-})}{vol_{d}(K)})^{k}
     (1- \frac{vol_{d}(K \cap
     H^{-})}{vol_{d}(K)})^{n-k}                         $$

$$ | \Bbb P_{K \cap H^{-}} \{(x_1,...,x_k)| x_t \notin [x_1,...,x_k] \}
   - \Bbb P_{B \cap \tilde H^{-}} \{(z_1,...,z_k)| z_t \notin
   [z_1,...,z_k] \} | + 2\epsilon  $$

\vskip 5mm

On the other hand, if we choose $t_{\epsilon}$ sufficiently small we
have for all $t \in (0,t_{\epsilon}]$ and all $k, 0 \leq k \leq
4dc^{d+2}$

$$ | \Bbb P_{K \cap H^{-}} \{(x_1,...,x_k)| x_t \notin [x_1,...,x_k] \}
   - \Bbb P_{B \cap \tilde H^{-}} \{(z_1,...,z_k)| z_t \notin
   [z_1,...,z_k] \} | <
   \epsilon
   \tag 18  $$

This finishes the proof. We establish now (18). For $k=0,...,d$ the
difference is trivially 0. Now we assume that $x=z$ and $N(x)=N(z)$. We
have that

$$  \Bbb P_{K \cap H^{-}} \{(x_1,...,x_k)| x_t \notin [\{x_1,...,x_k\}
    \cap B \cap \tilde H^{-}] \}=  \tag 19 $$

$$ \sum_{m=0}^{k} \binom{k}{m} (\frac{vol_d(K \cap H^{-} \cap B \cap
\tilde H^{-})}
   {vol_{d}(K \cap H^{-})})^{m} (1- \frac{vol_{d}(K \cap H^{-} \cap B
   \cap \tilde H^{-})}{vol_{d}(K \cap H^{-})})^{k-m}             $$

$$  \Bbb P_{K \cap H^{-} \cap B \cap \tilde H^{-}}
   \{(x_1,...,x_{m})| x_t \notin [x_1,...,x_{m}]
   \}                       $$

\vskip 5mm

If we choose $t_{\epsilon}$ small enough we have by Lemma 14 for all $t
\in (0,t_{\epsilon}]$ that

$$ 1- \frac{vol_{d}(K \cap H^{-} \cap B \cap \tilde H^{-})}
   {vol_{d}(K \cap H^{-})}                                           $$

is so small that we get by (19) and $k \leq 4dc^{d+2}$

$$ | \Bbb P_{K \cap H^{-}} \{(x_1,...,x_k)| x_t \notin [\{x_1,...,x_k\}
\cap B \cap \tilde H^{-}] \}
\leq                                               \tag 20 $$

$$  \Bbb P_{K \cap H^{-} \cap B \cap \tilde H^{-}}
   \{(x_1,...,x_{k})| x_t \notin [x_1,...,x_{k}] \}+
   e^{-16dc}           $$

Moreover, we have

$$  \Bbb P_{K \cap H^{-} \cap B \cap \tilde H^{-}}
   \{(x_1,...,x_{k})| x_t \notin [x_1,...,x_{k}] \}        $$

$$ \geq \Bbb P_{K \cap H^{-} \cap B \cap \tilde H^{-}}
   \{(x_1,...,x_{k})|\{x_1,...,x_n\} \subset K \cap H^{+}(x_{t},N(x))
   \}   \tag 21 $$

$$   \geq (1- c^{-d-1})^{4dc^{d+2}} \geq e^{-8dc}
$$

The last inequality holds because we have $1-\frac{1}{s} \geq e^{-2s}$
for $s \geq 2$. By (20) and (21) we get now

$$ \Bbb P_{K \cap H^{-}} \{(x_1,...,x_k)| x_t \notin [\{x_1,...,x_k\}
\cap B \cap \tilde H^{-}] \}  $$ $$ \leq (1+ \exp(-c)) \Bbb P_{K \cap
H^{-} \cap B \cap \tilde H^{-}}
   \{(x_1,...,x_{k})| x_t \notin [x_1,...,x_{k}] \}           $$

Therefore we get now

$$ \Bbb P_{K \cap H^{-}} \{(x_1,...,x_k)| x_t \notin [x_1,...,x_k
    ] \}
    \leq
    $$

$$ \Bbb P_{K \cap H^{-}} \{(x_1,...,x_k)| x_t \notin [\{x_1,...,x_k\}
    \cap B \cap \tilde H^{-}]
    \}                                            $$

$$ \leq (1+e^{-c}) \Bbb P_{K \cap H^{-} \cap B \cap \tilde H^{-}}
   \{(x_1,...,x_{k})| x_t \notin [x_1,...,x_{k}] \}           $$

$$ \leq (1+ e^{-c})( \frac{vol_{d}(K \cap H^{-})} {vol_d(K \cap H^{-}
\cap B \cap \tilde H^{-})})^{k}            $$

$$  \Bbb P_{K \cap H^{-}} \{(x_1,...,x_k)| x_t \notin [x_1,...,x_k]
    \text{\ and\ } \{x_1,...,x_k\} \subset B \cap \tilde H^{-}
    \}                    $$

$$ \leq (1+ e^{-c})(\frac{vol_{d}(K \cap H^{-})} {vol_d(K \cap H^{-}
\cap B \cap \tilde H^{-})})^{k}
    \Bbb P_{K \cap H^{-}} \{(x_1,...,x_k)| x_t \notin [x_1,...,x_k]
    \}      $$

Thus we get that

$$ |\Bbb P_{K \cap H^{-}} \{(x_1,...,x_k)| x_t \notin [x_1,...,x_k
    ] \} - \Bbb P_{K \cap H^{-} \cap B \cap \tilde H^{-}}
   \{(x_1,...,x_{k})| x_t \notin [x_1,...,x_{k}] \}| < \epsilon
   $$

if we choose c sufficiently big. We have the same inequality for $\Bbb
P_{B \cap \tilde H^{-}}$.  This implies (18).  \qed \enddemo \vskip 1cm

\proclaim{\smc Lemma 18} Let K be a convex body in $\Bbb R^d$ and $x
\in \partial K$. Suppose that $\partial K$ is twice differentiable at x
in the generalized sense. Then we have

$$(i) \lim_{t \to 0} \frac{<x,N(x)>}{<x_t,N(x_t)>} =1  $$

\leftline{(ii)}

$$ \lim_{t \to 0} \frac{t^{\frac{d-1}{d+1}}}{vol_{d-1}(P_{N(x_t)}
   ((-x_t+K) \cap (x_t-K)))} = {\kappa(x)}^{\frac{1}{d+1}}
   (\frac{2}{d+1})^ {\frac{d-1}{d+1}}
   vol_{d-1}(B_{2}^{d-1})^{-\frac{2}{d+1}}          $$

\endproclaim \vskip 1cm

\demo{\smc Proof} (i) The same arguments as in the proof of Lemma 5 are
applied. We just sketch the argument. Suppose (i) is not true. Then we
find a supporting hyperplane $H(x_t,N(x_t))$ so that $x_t$ is very
close to x but $N(x_t)$ is not close to $N(x)$. By the assumption we
have that all the points in the set $H(x_t,N(x_t)) \cap K$ do not
belong to the interior of $K_t$. On the other hand, the volume
$vol_{d}(K \cap H^{-}(x_t,N(x_t))$ is so big that we can single out a
point in $H(x_t,N(x_t)) \cap K$ that is in the interior of $K_t$.  \par

(ii) We consider the case $\kappa(x) > 0$. The case $\kappa(x)=0$ is
treated in an analogous way. By Lemma 14 K can be approximated by an
ellipsoid. We may assume it is a Euclidean sphere. By (i)
$<x_t,N(x_t)>$ is as close to $<x,N(x)>$ as we choose it to be for
small t. Altogether we have that

$$ vol_{d-1}(P_{N(x_t}
   ((-x_t+K) \cap (x_t-K)))          $$

is up to some error equal to

$$ vol_{d-1}(H(x_t,N(x)) \cap K)     $$

or

$$ vol_{d-1}(H(x_t,N(x)) \cap
B_{2}^{d}(x-{\kappa(x)}^{-\frac{1}{d-1}}N(x),
   {\kappa(x)}^{-\frac{1}{d-1}})  $$

It is left to apply Lemma 4.  \qed \enddemo \vskip 1cm

\demo{\smc Proof of Theorem 1} We may assume that $K_T$ coincides with
the origin. By Lemma 3 and 8 we have

$$ lim_{n \to \infty} \frac{vol_{d}(K)- \Bbb E(K,n)}{(\frac{1}
   {n})^{\frac{2}{d+1}}}=
   $$

$$ \lim_{n \to \infty} n^{\frac{2}{d+1}} \int_{\partial K}
\int_0^{\frac{T}{2}}
 \frac{\Bbb P \{(x_1,...., x_n) \mid x_t \notin [x_1,....,x_n] \} }
{vol_{d-1}(P_{N(x_t)}((-x_t+K) \cap (x_t-K)))} \frac {\Vert x_t
\Vert^d}{\Vert x \Vert ^d}  \frac {<x,N(x)>}{<x_t,N(x_t)>} dt d\mu $$

provided the limit exists. We apply now Lebesgue's convergence theorem
in order to change limit and integration over $\partial K$. The
hypothesis of Lebesgue's theorem is fulfilled because of Lemma 6 and
10. By Lemma 12 we get that the latter expression equals

$$  \int_{\partial K}\lim_{n \to \infty} n^{\frac{2}{d+1}}
\int_0^{\frac{\log n}{n}}
 \frac{\Bbb P \{(x_1,...., x_n) \mid x_t \notin [x_1,....,x_n] \} }
{vol_{d-1}(P_{N(x_t)}((-x_t+K) \cap (x_t-K)))} \frac {\Vert x_t
\Vert^d}{\Vert x \Vert ^d}  \frac {<x,N(x)>}{<x_t,N(x_t)>} dt d\mu $$

By Lemma 18 this expression equals

$$  \int_{\partial K}\lim_{n \to \infty} n^{\frac{2}{d+1}}
\int_0^{\frac{\log n}{n}}
 \frac{\Bbb P \{(x_1,...., x_n) \mid x_t \notin [x_1,....,x_n] \} }
{vol_{d-1}(P_{N(x_t)}((-x_t+K) \cap (x_t-K)))} dt d\mu $$

By Lemma 18 (ii) we get

$$ \int_{\partial K} {\kappa(x)}^{\frac{1}{d+1}} d\mu
   \lim_{n \to \infty} n^{\frac{2}{d+1}} \frac{
   (\frac{2}{d+1})^{\frac{d-1}{d+1}}}{
   {vol_{d-1}(B_2^{d-1})}^{\frac{2}{d+1}}} \int_{0}^{\frac{\log n}{n}}
   \frac{ \Bbb P \{(x_1,...,x_n)| x_t \notin [x_1,...,x_n]\}}
   {t^{\frac{d-1}{d+1}}}  dt                                 $$

By Lemma 17 we have for $x \in \partial K$ with $\kappa(x) >0$

$$ \lim_{n \to \infty} n^{\frac{2}{d+1}} \int_{0}^{\frac{\log n}{n}}
t^{-\frac{d-1}{d+1}}
   \Bbb P_{K} \{(x_1,...,x_n)| x_t \notin [x_1,...,x_n]\}
   dt                 $$

$$ =\lim_{n \to \infty} n^{\frac{2}{d+1}} \int_{0}^{\frac{\log n}{n}}
t^{-\frac{d-1}{d+1}}
   \Bbb P_{B} \{(z_1,...,z_n)| z_t \notin [z_1,...,z_n]\}
   dt                 $$

where B is a Euclidean ball whose volume is the same as that of K.  The
limit for B exists by Lemma 13. Thus we get

$$ lim_{n \to \infty} \frac{vol_{d}(K)- \Bbb E(K,n)}{(\frac{1}
   {n})^{\frac{2}{d+1}}}=
   $$

$$ \int_{\partial K} {\kappa(x)}^{\frac{1}{d+1}} d\mu
   \lim_{n \to \infty} n^{\frac{2}{d+1}}
   \frac{(\frac{2}{d+1})^{\frac{d-1}{d+1}}}{
   {vol_{d-1}(B_2^{d-1})}^{\frac{2}{d+1}}} \int_{0}^{\frac{\log n}{n}}
   \frac{ \Bbb P_{B} \{(z_1,...,z_n)| z_t \notin [z_1,...,z_n]\}}
    {t^{\frac{d-1}{d+1}}} dt                                        $$

Since this formula holds for all convex bodies it holds in particular
for the Euclidean ball. By Lemma 13 we determine the coefficient.  \qed
\enddemo

 \enddocument


\Refs \widestnumber\key{Schn-Wie}

\ref \key A \by A.D. Aleksandrov
   \paper Almost everywhere existence of the second differential of a
   convex function and some properties of convex surfaces connected
   with it \jour Uchenye Zapiski Leningrad Gos. Univ., Math. Ser.  \yr
   1939 \vol 6  \pages 3--35 \endref

\ref \key Ba \by V. Bangert
   \paper Analytische Eigenschaften konvexer Funktionen auf
   Riemannschen Mannigfaltigkeiten \jour J. Reine Angew. Math. \yr 1979
   \vol 307 \pages 309--324 \endref

\ref \key B \by I. B\'ar\'any \pages 81--92
   \paper Pandom polytopes in smooth convex bodies \yr 1992 \vol 39
   \jour Mathematika \endref

\ref \key B-L \by I. B\'ar\'any and D.G. Larman \pages 274--291
   \paper Convex bodies, economic cap covering, random polytopes \yr
   1988 \vol 35 \jour Mathematika  \endref

\ref \key F-R \by I. F\'ary and L. R\'edei \pages 205--220
   \paper Der zentralsymmetrische Kern und die zentralsymmetrische
    H\"ulle von konvexen K\"orpern \yr 1950 \vol 122 \jour Math. Ann.
   \endref

\ref \key G \by B. Gr\"unbaum
   \paper Measures of symmetry for convex sets \inbook Convexity,
   Proceedings of Symposion in Pure Mathematics \publ AMS \ed V.L. Klee
   \yr 1963 \pages 233--270 \endref

\ref \key K \by K. Kiener
   \paper Extremalit\"at von Ellipsoiden und die Faltungsungleichung
   von Sobolev \jour Arch. Math. \yr 1986 \vol 46 \pages 162--168
   \endref

\ref \key R-S1 \by A. R\'enyi and R. Sulanke
   \paper \"Uber die konvexe H\"ulle von n zuf\"allig gew\"ahlten
   Punkten \jour Zeitschr. Wahrschein. verwandte Geb.  \yr 1963 \vol 2
   \pages 75--84 \endref

\ref \key R-S2 \by A. R\'enyi and R. Sulanke
   \paper \"Uber die konvexe H\"ulle vo n zuf\"allig gew\"ahlten
   Punkten II \jour Zeitschr. Wahrschein. verwandte Geb.  \yr 1964 \vol
   3 \pages 138--147 \endref

\ref \key Schm1 \by M. Schmuckenschl\"ager
   \paper The distribution function of the convolution square of a
   convex symmetric body in $R^n$ \jour Israel J. Math. \yr 1992 \vol
   78 \pages 309--334 \endref

\ref \key Schm2 \by M. Schmuckenschl\"ager
   \paper Notes \endref

\ref \key Schn \by R. Schneider \pages 211--227
   \paper Random approximation of convex sets \yr 1988 \vol 151 \jour
   J. Microscopy \endref

\ref \key Schn-Wie \by R. Schneider and J.A. Wieacker \pages 69--73
   \paper Random polytopes in a convex body \yr 1980 \vol 52 \jour Z.
   Wahrscheinlichkeitstheorie verw. Gebiete \endref

\ref \key Sch\"u-W 1 \by C. Sch\"utt and E. Werner \pages 275--290
   \paper The convex floating body \yr 1990 \vol 66 \jour Math. Scand.
   \endref

\ref \key Sch\"u-W 2 \by C. Sch\"utt and E. Werner \pages 169--188
   \paper The convex floating body of almost polygonal bodies \yr 1992
   \vol 44 \jour Geometriae Dedicata \endref

\ref \key St \by S. Stein \pages 145--148
   \paper The symmetry function in a convex body \yr 1956 \vol 6 \jour
   Pacific J. Math. \endref

\ref \key Wie \by J.A. Wieacker
   \book Einige Probleme der polyedrischen Approximation \bookinfo
   Diplomarbeit \publaddr Frei burg im Breisgau \yr 1978 \endref


\endRefs
\bye